\documentclass  [11pt] {article}
\usepackage [margin=.8in]{geometry}
\usepackage[utf8]{inputenc}
\usepackage{amsmath,amsthm,amssymb,amsopn}
\usepackage{authblk}
\usepackage{lipsum}

\newcommand\blfootnote[1]{%
  \begingroup
  \renewcommand\thefootnote{}\footnote{#1}%
  \addtocounter{footnote}{-1}%
  \endgroup
}

\begin{document}

\title{Condensable models of set theory}
\author{Ali Enayat}
\maketitle

\begin{abstract}
\noindent A model $\mathcal{M}$ of ZF is said to be \textit{condensable }if $%
\mathcal{M}\cong \mathcal{M(\alpha )}\prec _{\mathbb{L}_{\mathcal{M}}}%
\mathcal{M}$ for some \textquotedblleft ordinal\textquotedblright\ $\alpha
\in \mathrm{Ord}^{\mathcal{M}}$, where $\mathcal{M(\alpha )}:=(\mathrm{V}(\alpha
),\in )^{\mathcal{M}}$ and $\mathbb{L}_{\mathcal{M}}$ is the set of formulae
of the infinitary logic $\mathbb{L}_{\infty ,\omega }$ that appear in the
well-founded part of $\mathcal{M}$. The work of Barwise and Schlipf in the
1970s revealed the fact that every countable recursively saturated model of
ZF is cofinally condensable (i.e., $\mathcal{M}\cong \mathcal{M(\alpha )}%
\prec _{\mathbb{L}_{\mathcal{M}}}\mathcal{M}$ for an unbounded collection of
$\alpha \in \mathrm{Ord}^{\mathcal{M}}$). Moreover, it can be readily shown
that any $\omega $-nonstandard condensable model of $\mathrm{ZF}$ is
recursively saturated. These considerations provide the context for the
following result that answers a question posed to the author by Paul
Kindvall Gorbow.\medskip

\noindent \textbf{Theorem A.} \textit{Assuming a modest set-theoretic
hypothesis, there is a countable model }$\mathcal{M}$ \textit{of }ZFC
\textit{that is both \textbf{definably well-founded} }(\textit{i.e., every first order definable element of
 }$\mathcal{M}$  \textit{\ is in the well-founded part of }$\mathcal{M)}$
\textit{and \textbf{cofinally condensable}.}\medskip

\noindent We also provide various equivalents of the notion of condensability, including the result below.\medskip

\noindent \textbf{Theorem B.} \textit{The following are equivalent for a
countable model} $\mathcal{M}$ \textit{of }$\mathrm{ZF}$:\medskip

\noindent \textbf{(a)} $\mathcal{M}$ \textit{is condensable}.\medskip

\noindent \textbf{(b) }$\mathcal{M}$ \textit{is cofinally condensable}%
.\medskip

\noindent \textbf{(c)} $\mathcal{M}$ \textit{is nonstandard and }$\mathcal{%
M(\alpha )}\prec _{\mathbb{L}_{\mathcal{M}}}\mathcal{M}$ \textit{for an
unbounded collection\ of }$%
\alpha \in \mathrm{Ord}^{\mathcal{M}}$.
\end{abstract}

\blfootnote {\textit{Acknowledgments}. I am grateful to Paul Kindvall Gorbow, Zachiri McKenzie, Jim Schmerl, and an anonymous referee for their feedback on earlier versions of this paper.}



\begin{center}
\textbf{1.~INTRODUCTION}\bigskip
\end{center}

By a classical theorem of Harvey Friedman \cite{Friedman}, every countable
nonstandard model $\mathcal{M}$ of $\mathrm{ZF}$ can be \textquotedblleft
shrunk\textquotedblright\ in the sense that $\mathcal{M}$ is isomorphic to a
proper rank-initial segment of itself\textit{. }Friedman's theorem was
refined by Jean-Pierre Ressayre \cite{Ressayre-monograph}, who constructed
proper rank-initial self-embeddings of models of set theory that pointwise
fix any prescribed rank-initial segment $\mathcal{M(\alpha )}:=(\mathrm{V}%
(\alpha ),\in )^{\mathcal{M}}$ of a model $\mathcal{M}$ of set theory
determined by an \textquotedblleft ordinal\textquotedblright\ $\alpha \in $ $%
\mathrm{Ord}^{\mathcal{M}}$. More recently, Paul Kindvall Gorbow \cite%
{Gorbow-Dissertation} extended Ressayre's work by carrying out a systematic
study of the structure of \textit{fixed point sets} of rank initial
self-embeddings of models of set theory; and Zachiri McKenzie and the author
\cite{Ali + Zach}, studied self-embeddings whose images are only required to
be $\in $-initial segments of the ambient models.\medskip

By a general result of Gorbow \cite[Theorem 7.2]{Gorbow-AML} if $\mathcal{M}$
is a countable nonstandard model of the fragment $\mathrm{KP}^{\mathcal{P}%
}+\Sigma _{1}^{\mathcal{P}}$-$\mathrm{Separation}$ of $\mathrm{ZF}$, then
there are continuum-many proper rank-initial segments $\mathcal{N}$ of $%
\mathcal{M}$ that are isomorphic to $\mathcal{M}$, which makes it clear that
there are continuum many such rank-initial segments $\mathcal{N}$ that are
not of the form $\mathcal{M(\alpha )}$ for any \textquotedblleft
ordinal\textquotedblright\ $\alpha \in \mathrm{Ord}^{\mathcal{M}}$%
(equivalently: $\mathrm{Ord}^{\mathcal{M}}\backslash \mathrm{Ord}^{\mathcal{N%
}}$ has no least element). It is also known that every consistent extension
of ZF has a countable model $\mathcal{M}$ that is not isomorphic to any
initial segment of itself that is of the form $\mathcal{M(\alpha )}$ for any
$\alpha \in \mathrm{Ord}^{\mathcal{M}}$.\footnote{%
This follows from the following facts: (1) Every countable extension of ZF
has a Paris model, i.e., a model every ordinal of which is pointwise
definable; (2) No Paris model of ZF can be rank extended to a model of
Kripke-Platek set theory that has a first new ordinal; see Theorems 1.2 and
3.11 of \cite{Enayat-DO}.} On the other hand, the pioneering work of Barwise
and Schlipf on recursively saturated models in the 1970s revealed a wealth
of countable models $\mathcal{M}$ of set theory such that $\mathcal{M\cong
M(\alpha )}\prec \mathcal{M}$ for many $\alpha \in \mathrm{Ord}^{\mathcal{M}%
},$ as indicated by the following theorem.\medskip

\noindent \textbf{1.1.}~\textbf{Theorem.}~(Barwise and Schlipf) \textit{Let }%
$\mathcal{M}$\textit{\ be a countable recursively saturated model of }$%
\mathit{\mathrm{ZF.}}$ \textit{Then there is some} $\alpha \in \mathrm{Ord}^{%
\mathcal{M}}$ \textit{such that} $\mathcal{M}\cong \mathcal{M(\alpha )}\prec
\mathcal{M}$. \textit{Moreover, the collection of such} $\alpha $\textit{s
is unbounded in }$\mathrm{Ord}^{\mathcal{M}}.$\footnote{%
\noindent Theorem 1.1 can be readily derived from the following two key
results of Barwise and Schlipf:
\par
\noindent (a) Every resplendent model of ZF is cofinally condensable \cite[%
Corollary 3.3]{Schlipf-JSL}.
\par
\noindent (b) Every countable recursively saturated model is resplendent
\cite{Barwise-Schlipf-JSL}.
\par
It is worth pointing out that the assumption of countability in Theorem 1.1
cannot be dropped since it is well-known that every consistent extension of
ZF has an $\aleph _{1}$-like recursively saturated model (as elaborated in
Remark 4.7). On the other hand, in light of the resplendence property of
special models (attributed to Chang and Moschovakis in \cite[Example 2.3]%
{Barwise-Schlipf-JSL}), every saturated model of ZF is cofinally condensable
(it is well-known that ZFC proves that a saturated model of ZF of
cardinality $\kappa $ exists iff ZF is consistent and $\kappa $ is an
uncountable cardinal such that $\kappa ^{<\kappa }=\kappa ).$}\medskip

\noindent Motivated by the above theorem, we say that a model $\mathcal{M}$
of ZF is \textit{condensable }if there is some $\alpha \in \mathrm{Ord}^{%
\mathcal{M}}$ such that $\mathcal{M}\cong \mathcal{M(\alpha )}\prec _{%
\mathbb{L}_{\mathcal{M}}}\mathcal{M},$ where $\mathcal{M(\alpha )}:=(\mathrm{%
V}(\alpha ),\in )^{\mathcal{M}}$ and $\mathbb{L}_{\mathcal{M}}$ consists of
formulae of the infinitary logic $\mathbb{L}_{\infty ,\omega }$ that appear
in the well-founded part of $\mathcal{M}$. Note that for an $\omega $%
-nonstandard model $\mathcal{M}$ of ZF, $\mathbb{L}_{\mathcal{M}}$ is just
the collection of (finitary) first order formulae, so the condition $%
\mathcal{M(\alpha )}\prec _{\mathbb{L}_{\mathcal{M}}}\mathcal{M}$ is
equivalent to $\mathcal{M(\alpha )}\prec \mathcal{M}$ for $\omega $%
-nonstandard models $\mathcal{M}$ (and in particular for recursively
saturated models $\mathcal{M}$). Theorem 1.1 provides the background for a
question posed in an e-mail to the author (June 14, 2019).\medskip

\noindent \textbf{1.2.~Question }(Gorbow).~\textit{Is there an }$\omega $%
\textit{-standard model }$\mathcal{M}$\textit{\ of }$\mathrm{ZF}$ \textit{%
such that the collection of }$\alpha \in \mathrm{Ord}^{\mathcal{M}}$ \textit{%
satisfying }$\mathcal{M}\cong \mathcal{M(\alpha )}\prec \mathcal{M}$ \textit{%
is unbounded in }$\mathrm{Ord}^{\mathcal{M}}$?\medskip

\noindent In this article we establish Theorem A of the abstract (in Section
3) and a strengthening of Theorem B of the abstract (in Section 4). Theorem
A yields a (strong) positive answer to Gorbow's question. Theorem B, on the
other hand, provides various characterizations of condensability over the
family of countable models of \textrm{ZF}. \bigskip

\begin{center}
\textbf{2.~PRELIMINARIES}\bigskip
\end{center}

In this section we collect the basic definitions, notations, conventions,
and results that will be used in the statements and proofs of our main
results in Sections 3 and 4.\medskip

\noindent \textbf{2.1.~Definition.}~(Models, languages, and theories) Models
will be represented by calligraphic fonts ($\mathcal{M}$, $\mathcal{N}$,
etc.) and their universes will be represented by the corresponding roman
fonts ($M$, $N$, etc.). In the definitions below, $\mathcal{M}$ is a model
of ZF and $\in ^{\mathcal{M}}$ is the membership relation of $\mathcal{M}$.
\medskip

\noindent \textbf{(a) }$\mathrm{Ord}^{\mathcal{M}}$ is the class of
\textquotedblleft ordinals\textquotedblright\ of $\mathcal{M}$, i.e., $%
\mathrm{Ord}^{\mathcal{M}}:=\left\{ m\in M:\mathcal{M}\models \mathrm{Ord}%
(m)\right\} ,$ where $\mathrm{Ord}(x)$ expresses \textquotedblleft $x$ is
transitive and is well-ordered by $\in $\textquotedblright . More generally,
given a class $\mathrm{D}$ whose defining formula is $\delta (x)$, $\mathrm{D%
}^{\mathcal{M}}:=\left\{ m\in M:\mathcal{M}\models \delta (m)\right\} .$%
\medskip

\noindent \textbf{(b) }$\mathcal{M}$ is \textit{nonstandard} if $\in ^{%
\mathcal{M}}$ is ill-founded (equivalently: if $(\mathrm{Ord},\in )^{%
\mathcal{M}}$ is ill-founded). $\mathcal{M}$ is $\omega $-\textit{nonstandard%
} if ($\omega ,\in )^{\mathcal{M}}$ is ill-founded.\medskip

\noindent \textbf{(c)} For $c\in M$, $\mathrm{Ext}_{\mathcal{M}}(c)$ is the $%
\mathcal{M}$-extension of $c$, i.e., $\mathrm{Ext}_{\mathcal{M}}(c):=\{m\in
M:m\in ^{\mathcal{M}}c\}.$ We say that a subset $X$ of $M$ \textit{is coded
in} $\mathcal{M}$ if there is some $c\in M$ such that $\mathrm{Ext}_{%
\mathcal{M}}(c)=X.$ For $A\subseteq M$, $\mathrm{Cod}_{A}\mathrm{(}\mathcal{%
M)}$ will denote the collection of sets of the form $A\cap \mathrm{Ext}_{%
\mathcal{M}}(c)$, where $c\in M$.\medskip

\noindent \textbf{(d)} Given a model $\mathcal{N}$ of ZF, we write $\mathcal{%
M}\subseteq _{\mathrm{rank}}\mathcal{N}$ (read as: $\mathcal{M}$ \textit{is
rank-extended by} $\mathcal{M}$, or $\mathcal{M}$ \textit{is a rank-initial
segment of} $\mathcal{N}$), to indicate that $\mathcal{M}$ is a submodel of $%
\mathcal{N}$ such that $\rho ^{\mathcal{M}}(x)>$ $\rho ^{\mathcal{M}}(m)$
whenever $x\in N\backslash M$ and $m\in M,$ where $\rho (x)$ is the usual
ordinal-valued rank function of set theory.\medskip

\noindent \textbf{(e)} The \textit{well-founded part} of $\mathcal{M}$,
denoted $\mathrm{WF}(\mathcal{M})$, consists of all elements $m$ of $%
\mathcal{M}$ such that there is no infinite sequence $\left\langle
a_{n}:n<\omega \right\rangle $ with $m=a_{0}$ and $a_{n+1}\in ^{\mathcal{M}%
}a_{n}$ for all $n\in \omega .$ Given $m\in M,$ we say that $m$ \textit{is a
nonstandard element of} $\mathcal{M}$\ if $m\notin \mathrm{WF}(\mathrm{M}).$
We denote the submodel of $\mathcal{M}$ whose universe is $\mathrm{WF}(%
\mathcal{M})$ by $\mathcal{WF}(\mathcal{M}).$ It is well-known that if $%
\mathcal{M}$ is a model of $\mathrm{ZF}$, then $\mathcal{WF}(\mathcal{M}%
)\subseteq _{\mathrm{rank}}\mathcal{M}$, and $\mathcal{WF}(\mathcal{M})$
satisfies KP (Kripke-Platek set theory) \cite[Chapter II, Theorem 8.4]%
{Barwise-book}. Also note that if $\mathcal{M}\subseteq _{\mathrm{rank}}%
\mathcal{N}$, where $\mathcal{M}$ is nonstandard, then $\mathcal{M}$ and $\mathcal{N}$ share the same
well-founded part.

\begin{itemize}
\item It is important to bear in mind that we will identify\textit{\ }$%
\mathcal{WF}(\mathcal{M})$ with its transitive collapse.
\end{itemize}

\noindent \textbf{(f)} $\mathrm{o}(\mathcal{M})$ (read as: \textit{the
ordinal of }$\mathcal{M}$) is the supremum of all ordinals that appear in
the well-founded part of $\mathcal{M}$.\medskip

\noindent \textbf{(g)} Let $\mathcal{L}_{\mathrm{set}}$ be the usual
vocabulary $\{=,\in \}$ of set theory. In this paper we use $\mathbb{L}%
_{\infty ,\omega }$ to denote the language using the vocabulary $\mathcal{L}%
_{\mathrm{set}}$ that allows conjunctions and disjunctions of \textit{sets}
(but not proper classes) of formulae, subject to the restriction that such
infinitary formulae have at most finitely many free variables. Given a set $%
\Psi $ of formulae, we denote such conjunctions and disjunctions
respectively as $\bigwedge \Psi $ and $\bigvee \Psi $.

\begin{itemize}
\item In the interest of efficiency, we will treat disjunctions and
universal quantification as defined notions.
\end{itemize}

\noindent \textbf{(h)} $\mathbb{L}_{\delta ,\omega }$ is the sublanguage of $%
\mathbb{L}_{\infty ,\omega }$ that allows conjunctions and disjunctions of
sets of formulae of cardinality \textit{less than} $\delta .$ Note that $%
\mathbb{L}_{\omega ,\omega }$ is none other than the usual first order
language of set theory, and that in general the language $\mathbb{L}_{\delta
,\omega }$ only uses finite strings of quantifiers (as indicated by the $%
\omega $ in the subscript). A full treatment can be found in \cite[Chapter
III]{Barwise-book}.$\medskip $

\noindent \textbf{(i) }Given $\mathbb{L}\subseteq \mathbb{L}_{\infty ,\omega
}$, and $\mathcal{L}_{\mathrm{set}}$-structures $\mathcal{N}_{1}$ and $%
\mathcal{N}_{2}$, we write $\mathcal{N}_{1}\prec _{\mathbb{L}}\mathcal{N}%
_{2} $ to indicate that $\mathcal{N}_{1}$ is a submodel of $\mathcal{N}_{2}$
and for every $\varphi (x_{1},\cdot \cdot \cdot ,x_{n})\in \mathbb{L}$, and
any $n$-tuple $(a_{1},\cdot \cdot \cdot ,a_{n})$ from $N_{1}$, we have:

\begin{center}
$\mathcal{N}_{1}\models \varphi (a_{1},\cdot \cdot \cdot ,a_{n})$ iff $%
\mathcal{N}_{2}\models \varphi (a_{1},\cdot \cdot \cdot ,a_{n})$.
\end{center}

\noindent \textbf{(j) }$\mathbb{L}_{\mathcal{M}}=\mathbb{L}_{\infty ,\omega
}\cap \mathrm{WF}(\mathcal{M)}$. Note that if $M$ is countable, \textbf{\ }$%
\mathbb{L}_{\mathcal{M}}=\mathbb{L}_{\omega _{1},\omega }\cap \mathrm{WF}(%
\mathcal{M)}$.\medskip

\noindent \textbf{(k)} Given $\mathbb{L}\subseteq \mathbb{L}_{\infty ,\omega
}$, $\mathrm{Th}_{\mathbb{L}}(\mathcal{M)}$ is the set of sentences (closed
formulae) of $\mathbb{L}$ that hold in $\mathcal{M}$, and $\mathrm{ZF}(%
\mathbb{L})$ is the natural extension of ZF in which the scheme $\mathrm{Sep}
$ of separation and $\mathrm{Coll}$ of collection are extended to the
schemes $\mathrm{Sep}(\mathbb{L})$ and $\mathrm{Coll}(\mathbb{L})$ to allow
formulae in $\mathbb{L}$ to be used for \textquotedblleft
separating\textquotedblright\ and \textquotedblleft
collecting\textquotedblright\ (respectively).\medskip

\noindent \textbf{(l)} For $\varphi \in \mathbb{L}_{\infty ,\omega }$, the
\textit{depth} of $\varphi ,$ denoted $\mathrm{Depth}(\varphi )$, is the
ordinal defined recursively by the following clauses:\medskip

\noindent (1) $\mathrm{Depth}(\varphi )=0$, if $\varphi $ is an atomic
formula.

\noindent (2) $\mathrm{Depth}(\varphi )=\mathrm{Depth}(\psi )+1,$ if $%
\varphi =\lnot \psi .$

\noindent (3) $\mathrm{Depth}(\varphi )=\mathrm{Depth}(\psi )+1,$ if $%
\varphi =\exists x\ \psi .$

\noindent (4) $\mathrm{Depth}(\varphi )=\sup \{\mathrm{Depth}(\psi )+1:\psi
\in \Psi \}$, if $\varphi =\bigwedge \Psi .$

\begin{itemize}
\item Within KP, one can code each formula $\varphi \in \mathbb{L}_{\infty
,\omega }$ with a set $\ulcorner \varphi \urcorner $ as in Chapter 3 of \cite%
{Barwise-book}, but in the interest of better readability we will often
identify a formula with its code. This coding allows us to construe
statements such as $\varphi \in \mathbb{L}_{\infty ,\omega }$ and $\mathrm{%
Depth}(\varphi )=\alpha $ as statements in the first order language of set
theory. It is easy to see that the collection $\mathrm{D(}\alpha )$ of
(codes of) $\mathbb{L}_{\infty ,\omega }$-formulae whose depth is less than $%
\alpha $ forms a set in ZF for all ordinals $\alpha $ since a simple
induction shows that for a sufficiently large $k\in \omega $, $\mathrm{D(}%
\alpha )\subseteq \mathrm{V(}\omega +k\alpha )$ for each ordinal $\alpha $.
This makes it clear that $\mathbb{L}_{\mathcal{M}}=\bigcup\limits_{\alpha
\in \mathrm{o}(\mathcal{M})}\mathrm{D}^{\mathcal{M}}\mathrm{(}\alpha ).$
\end{itemize}

\noindent \textbf{2.2.~Definition.~ }Suppose $\mathcal{M}$ is a model of $%
\mathrm{ZF}$, and $S\subseteq M$. \medskip

\noindent \textbf{(a)} $S$ is \textit{separative} (over $\mathcal{M}$) if $(%
\mathcal{M},S)$ satisfies the separation scheme $\mathrm{Sep(S)}$ in the
extended language that includes a fresh predicate \textrm{S} (interpreted by
$S$). \medskip

\noindent \textbf{(b)} $S$ is \textit{collective} (over $\mathcal{M}$) if $(%
\mathcal{M},S)$ satisfies the collection scheme $\mathrm{Coll(S)}$ in the
extended language that includes a fresh predicate \textrm{S} (interpreted by
$S$).\medskip

\noindent \textbf{(c)} $S$ is \textit{amenable} (over $\mathcal{M}$) if $S$
is both separative and collective$.$ In other words, $S$ is amenable if $%
\left( \mathcal{M},S\right) $ satisfies the replacement scheme $\mathrm{Repl}%
(S)$ in the extended language that includes a fresh predicate \textrm{S}
(interpreted by $S$). \medskip

\noindent \textbf{(d)} For $\alpha \in \mathrm{Ord}^{\mathcal{M}}$, $S$ is
an $\alpha $\textit{-satisfaction class} (over $\mathcal{M}$) if $S$
correctly decides the truth of atomic sentences, and $S$ satisfies Tarski's
compositional clauses of a truth predicate for \textrm{D}$^{\mathcal{M}}$%
\textrm{(}$\alpha )$-sentences (see below for the precise definition). $S$
is an $\infty $\textit{-satisfaction class} over $\mathcal{M}$, if $S$ is an
$\alpha $\textit{-}satisfaction class over $\mathcal{M}$ for every $\alpha
\in \mathrm{Ord}^{\mathcal{M}}$. \medskip

\noindent We elaborate the meaning of (d) above. Reasoning within \textrm{KP}%
, for each object $a$ in the universe of sets, let $c_{a\text{ }}$ be a
constant symbol denoting $a$ (where the map $a\mapsto c_{a}$ is $\Delta
_{1}),$ and let $\mathrm{Sent}^{+}(\alpha ,x)$ be the set-theoretic formula
(with an ordinal parameter $\alpha $ and the free variable $x)$ that defines
the proper class of sentences of the form $\varphi \left( c_{a_{1}},\cdot
\cdot \cdot ,c_{a_{n}}\right) $, where $\varphi (x_{1},\cdot \cdot \cdot
,x_{n})\in \mathrm{D}$\textrm{(}$\alpha )$ (the superscript $+$ on $\mathrm{%
Sent}^{+}(\alpha ,x)$ indicates that $x$ is a sentence in the language
augmented with the indicated proper class of constant symbols). Then $S$ is
an $\alpha $-satisfaction class over $\mathcal{M}$ if $\left( \mathcal{M}%
,S\right) \models \mathrm{Sat}(\alpha ,\mathrm{S})$, where $\mathrm{Sat}%
(\alpha ,\mathrm{S})$ is the (universal generalization of) the conjunction
of the axioms $(I)$ through $(IV)$ below.\medskip

\noindent $(I)\ \ \left( \left( \mathrm{S}\left( \ulcorner
c_{a}=c_{b}\urcorner \right) \leftrightarrow a=b\right) \wedge \left(
\mathrm{S}\left( \ulcorner c_{a}\in c_{b}\urcorner \right) \leftrightarrow
a\in b\right) \right) .$\medskip

\noindent $(II)\ \ \left( \mathrm{Sent}^{+}(\alpha ,\varphi )\wedge \left(
\varphi =\lnot \psi \right) \right) \rightarrow \left( \mathrm{S}(\varphi
)\leftrightarrow \lnot \mathrm{S}\mathsf{(}\psi \mathsf{)}\right) \mathsf{.}$%
\medskip

\noindent $(III)$ $\ \left( \mathrm{Sent}^{+}(\alpha ,\varphi )\wedge \left(
\varphi =\bigwedge \Psi \right) \right) \rightarrow \left( \mathrm{S}%
(\varphi )\leftrightarrow \forall \psi \in \Psi \ \mathrm{S}\mathsf{(}\psi
\mathsf{)}\right) \mathsf{.}$\medskip

\noindent $(IV)$ $\ \left( \mathrm{Sent}^{+}(\alpha ,\varphi )\wedge \left(
\varphi =\exists x\ \psi (x)\right) \right) \rightarrow \left( \mathrm{S}%
(\varphi )\leftrightarrow \exists x\ \mathrm{S}\mathsf{(\psi (}c_{x}\mathsf{%
))}\right) .$\medskip

\noindent \textbf{(e)} For $\alpha <\mathrm{o}(\mathcal{M})$, $S$ is \textit{%
the} $\alpha $-satisfaction class over $\mathcal{M},$ if $S$ is the usual
Tarskian satisfaction class for formulae in $\mathbb{L}_{\mathcal{M}}$ of
depth less than $\alpha ,$ i.e., the unique $\alpha $-satisfaction class $S$
over $\mathcal{M}$ such that $S$ satisfies:\medskip

\noindent $(V)$ $\ \forall x\left( \mathrm{S}(x)\rightarrow \mathrm{Sent}%
^{+}(\alpha ,x)\right) .$\medskip

\begin{itemize}
\item In the interest of a lighter notation we will often confuse constant
symbols in formulae with their denotations, e.g., we will write $\varphi (a)$
instead of $\varphi (c_{a}).$
\end{itemize}

The following proposition is immediately derivable from the definitions
involved.\medskip

\noindent \textbf{2.3.~Proposition.}~\textit{If} $S$ \textit{is an} $\alpha $%
\textit{-satisfaction class} \textit{over} $\mathcal{M}\models \mathrm{KP}$
\textit{for some nonstandard ordinal} $\alpha $ \textit{of} $\mathcal{M}$,
\textit{then for all n-ary formula }$\varphi (x_{1},\cdot \cdot \cdot
,x_{n}) $ \textit{of} $\mathbb{L}_{\mathcal{M}}$ \textit{and all }$n$-%
\textit{tuples} $(a_{1},\cdot \cdot \cdot ,a_{n})$ \textit{from} $M$, we
have:

\begin{center}
$\mathcal{M}\models \varphi (a_{1},\cdot \cdot \cdot ,a_{n})$ iff $\varphi
(a_{1},\cdot \cdot \cdot ,a_{n})\in S.$
\end{center}

\noindent \textit{In particular, for all sentences} $\varphi $ \textit{of} $%
\mathbb{L}_{\mathcal{M}},$ $\varphi \in S$ \textit{iff} $\varphi \in \mathrm{%
Th}_{\mathbb{L}_{\mathcal{M}}}(\mathcal{M)}.$\medskip

\noindent \textbf{2.4.~Remark.}~Reasoning within ZF, given any limit ordinal
$\gamma ,$ $\left( \mathrm{V}(\gamma ),\in \right) $ carries a separative $%
\gamma $-satisfaction class $S$ since we can take $S$ to be the Tarskian
satisfaction class on $\left( \mathrm{V}(\gamma ),\in \right) $ for formulae
of depth less than $\gamma $. More specifically, the Tarski recursive
construction/definition of truth works equally well in this more general
context of infinitary languages since $\left( \mathrm{V}(\gamma ),\in
\right) $ forms a set. Observe that $\left( \mathrm{V}(\gamma ),\in
,S\right) \models \mathrm{Sep(S)}$ comes \textquotedblleft for
free\textquotedblright\ since for any $X\subseteq \mathrm{V}(\gamma )$ the
expansion $\left( \mathrm{V}(\gamma ),\in ,X\right) $ satisfies the scheme
of separation in the extended language. However, $S$ need not be collective,
for example if the collection of first order definable elements of $\left(
\mathrm{V}(\gamma ),\in \right) $ is cofinal in $\gamma $.\medskip

\noindent \textbf{2.5.~Proposition.}~\textit{If}$\mathcal{M}\models \mathrm{%
KP}$, \textit{then for each }$\alpha \in \mathrm{o}(\mathcal{M})$ \textit{%
there is a formula} $\mathrm{Sat}_{\alpha }(x)\in \mathbb{L}_{\mathcal{M}}$
\textit{such that }$\mathrm{Sat}_{\alpha }^{\mathcal{M}}(x)$ \textit{is the }%
$\alpha $-\textit{satisfaction class over} $\mathcal{M}.$ \medskip

\noindent \textbf{Proof.} The desired formula $\mathrm{Sat}_{\alpha }(x)$ is
defined by the following recursion. A routine induction on $\alpha $ shows
that $\mathrm{Sat}_{\alpha }(x)$ has the desired properties. One can also
verify that $\mathrm{Sat}_{\alpha }(x)\in \mathrm{WF}(\mathcal{M})$ for $%
\alpha \in \mathrm{o}(\mathcal{M})$ by observing that the recursion defining
$\mathrm{Sat}_{\alpha }(x)$ is a so-called $\Sigma _{1}$-recursion (recall
that $\mathcal{WF}(\mathcal{M})$\ satisfies $\mathrm{KP}$, and $\mathrm{KP}$
can handle constructions by $\Sigma _{1}$-recursion).

\begin{itemize}
\item $\mathrm{Sat}_{1}(x):=\exists y\exists z\left[ \left( \left(
x=\ulcorner c_{y}=c_{z}\urcorner \right) \wedge (y=z)\right) \vee \left(
\left( x=\ulcorner c_{y}\in c_{z}\urcorner \right) \wedge (y\in z)\right) %
\right] .$

\item For $\alpha >1,$ $\mathrm{Sat}_{\alpha }(x):=\left[ \left( \mathrm{%
Depth}(x)=0\right) \wedge \mathrm{Sat}_{1}(x)\right] \lor$

\begin{equation*}
\bigvee\limits_{0<\beta <\alpha }\left( \mathrm{Depth}(x)=\beta \wedge
\lbrack \mathrm{Neg}_{\beta }(x)\vee \mathrm{Exist}_{\beta }(x)\vee \mathrm{%
Conj}_{\beta }(x)]\right),
\end{equation*}
where:
\end{itemize}

\begin{center}
$\mathrm{Neg}_{\beta }(x):=\exists y\left [\left( x=\ulcorner \lnot
y\urcorner \right) \wedge \lnot \mathrm{Sat}_{\beta }(y)\right],$\medskip

$\mathrm{Exist}_{\beta }(x):=\exists y\ \left[\exists v\left( x=\ulcorner
\exists v\ y(v)\urcorner \right) \wedge \exists v\ \mathrm{Sat}_{\beta
}(y(c_{v})\right]$, and

\medskip

$\mathrm{Conj}_{\beta }(x):=\exists y \left[ \left( x=\ulcorner \bigwedge
y\urcorner )\right) \wedge \left( \forall z\in y\ \mathrm{Sat}_{\beta
}(z)\right) \right] .$
\end{center}

\hfill $\square $\medskip

The following proposition will be called upon in the proofs of Theorems A
and B.\medskip

\noindent \textbf{2.6.} \textbf{Proposition.} (Overspill) \textit{Suppose} $%
\mathcal{M}$ \textit{is a nonstandard model of} $\mathrm{ZF}$, \textit{and} $%
S\subseteq M$ \textit{such that} $S$ \textit{is separative over} $\mathcal{M}
$. \textit{Assume furthermore that there is a first order formula} $\theta
(x,\overline{y})$ \textit{in the language} $\left\{ \in ,\mathrm{S}\right\} $
\textit{and some sequence of parameters} $\overline{m}\in M$ \textit{such
that }$(\mathcal{M},S)\models \theta (\alpha ,\overline{m})$ \textit{for
every} $\alpha \in \mathrm{o}(\mathcal{M}).$ \textit{Then there is a
nonstandard} $\gamma \in \mathrm{Ord}^{\mathcal{M}}$ \textit{such that} $(%
\mathcal{M},S)\models \theta (\gamma ,\overline{m}).\medskip $

\noindent \textbf{Proof.} Suppose not, and let $A:=\mathrm{WF}(\mathcal{M}%
)\cap \mathrm{Ord}^{\mathcal{M}}.$ Then $A=\left\{ x\in M:\left( \mathcal{M}%
,S\right) \models \theta (x,\overline{m})\wedge \mathrm{Ord}(x)\right\} $.
Since $A$ is a bounded subset of $\mathrm{Ord}^{\mathcal{M}}$, $A$ is coded
in $\mathcal{M}$ by $\mathrm{Sep(S)}$, and therefore has a supremum $\beta $
in $\mathcal{M}$. This is a contradiction since $\left( \beta ,\in \right) ^{%
\mathcal{M}}$ is well-founded, and yet $\beta \notin A$ since $A$ has no
last element.\hfill $\square $\medskip

The following general versions of the elementary chains theorem of model
theory, and of the reflection theorem of set theory will be called upon in
the proof of Theorem B. The proofs of Proposition 2.7 is obtained by a
routine adaptations of the usual proofs of the $\mathbb{L}_{\omega ,\omega }$%
-version of the elementary chain theorem (e.g., as in \cite[Theorem 3.1.9]%
{Chang-Keisler}). \medskip

\noindent \textbf{2.7.~Proposition.}~(Elementary Chains) \textit{Suppose }$%
\mathbb{L}\subseteq \mathbb{L}_{\infty ,\omega }$ \textit{where }$\mathbb{L}$
\textit{is closed under subformulae;} $(I,\vartriangleleft )$ \textit{is a
linear order; }$\left\langle \mathcal{M}_{i}:i\in I\right\rangle $ \textit{%
is an} $\mathbb{L}$-\textit{elementary chain of structures} (\textit{i.e. ,}
$\mathcal{M}_{i}\prec _{\mathbb{L}}\mathcal{M}_{j}$ \textit{whenever} $%
i\vartriangleleft j$); \textit{and} $\mathcal{M}=\bigcup\limits_{i\in I}%
\mathcal{M}_{i}.$\textit{\ Then }$\mathcal{M}_{i}\prec _{\mathbb{L}}\mathcal{%
M}$ \textit{for each} $i\in I.$\medskip

\noindent \textbf{2.8.~Proposition.}~(Reflection)\textbf{\ }\textit{Suppose }%
$\mathcal{M}\models \mathrm{ZF}(\mathbb{L}_{\mathcal{M}}),$ \textit{and for
each} $\varphi \in \mathbb{L}_{\mathcal{M}}$ \textit{where} $\varphi $
\textit{is} $n$\textit{-ary, let }$\mathrm{Ref}_{\varphi }(\gamma )$ \textit{%
be the }$\mathbb{L}_{\mathcal{M}}$-\textit{formula}:

\begin{center}
$\forall x_{1}\in \mathrm{V}(\gamma )\cdot \cdot \cdot \forall x_{n}\in
\mathrm{V}(\gamma )\ \left[ \varphi \left( x_{1},\cdot \cdot \cdot
,x_{n}\right) \longleftrightarrow \varphi ^{\mathrm{V}(\gamma )}\left(
x_{1},\cdot \cdot \cdot ,x_{n}\right) \right]. $
\end{center}

\noindent \textit{Then for any} $\alpha \in \mathrm{o}(\mathcal{M})$ \textit{%
there are arbitrarily large} $\gamma \in \mathrm{Ord}^{\mathcal{M}}$ \textit{%
such that }$\mathcal{M}(\gamma )$ \textit{reflects all formulae in} \textrm{D%
}$^{\mathcal{M}}$($\alpha ),$ \textit{i.e}., $\mathcal{M}\models \mathrm{Ref}%
_{\varphi }(\gamma )$ \textit{for all }$\mathbb{L}_{\mathcal{M}}$\textit{%
-formulae} $\varphi $ \textit{of depth less than} $\alpha .\medskip $

\noindent \textbf{Proof.} We will take advantage of Proposition 2.5 to
derive Proposition 2.8 from the usual finitary formulation of Reflection
Theorem\footnote{%
The Reflection Theorem is often formulated as a theorem scheme of ZF (e.g.,
as in \cite{Jechbook}), but it is well-known that the proof strategy of the
Reflection Theorem applies equally well to the extension $\mathrm{ZF}(%
\mathcal{L})$ of \textrm{ZF} (where $\mathcal{L}$ extends the language of
set theory) in which the schemes of separation and collection are extended
to $\mathcal{L}$-formulae.}. Given $\alpha \in \mathrm{o}(\mathcal{M}),$ by
Proposition 2.5, there is a formula $\mathrm{Sat}_{\alpha }(x)\in \mathbb{L}%
_{\mathcal{M}}$ such that\textit{\ }$\mathrm{Sat}_{\alpha }^{\mathcal{M}}(x)$
is an $\alpha $-satisfaction class over $\mathcal{M}$. Let $S:=\mathrm{Sat}%
_{\alpha }^{\mathcal{M}}(x)$; note that since $\mathcal{M}$ is assumed to
satisfy $\mathrm{ZF}(\mathbb{L}_{\mathcal{M}})$, $S$ is amenable over $%
\mathcal{M}$. Let $\mathrm{Sat}(\alpha ,\mathrm{S})$ be the finitary formula
that express \textquotedblleft \textrm{S} is an $\alpha $-satisfaction
class\textquotedblright\ (as in part (d) of Definition 2.2). By the
amenability of $S$ over $\mathcal{M}$, we may invoke the usual finitary
Reflection Theorem to obtain arbitrarily large $\gamma \in \mathrm{Ord}^{M}$
such that $\left( \mathcal{M}(\gamma ),S\cap M(\gamma )\right) $ reflects
the formula $\psi (x,\mathrm{S}):=\mathrm{Sat}(\alpha ,\mathrm{S})\wedge
\mathrm{S}(x),$ i.e.,$\medskip $

\begin{center}
$\left( \mathcal{M},S\right) \models \forall x\in \mathrm{V}(\gamma )\ \left[
\psi (x,\mathrm{S})\longleftrightarrow \psi ^{\mathrm{V}(\gamma )}(x,\mathrm{%
S})\right] .$
\end{center}

\noindent It is now easy to verify, using Proposition 2.3, that this choice
of $\gamma $ satisfies the conclusion of Proposition 2.8. \hfill $\square $%
\bigskip

\begin{center}
\textbf{3.~PROOF OF THEOREM A}\bigskip
\end{center}

The proof of Theorem A is based on some preliminary results. The reader may
wish to skip the proofs of these results in the first reading to be able to
better see the overall structure of the proof of Theorem A. \medskip

\noindent \textbf{3.1.~Lemma.}~\textit{Suppose }$\mathcal{M}$ \textit{and} $%
\mathcal{N}$ \textit{are countable nonstandard models of} $\mathrm{ZF}$
\textit{with the same well-founded part} $W$, \textit{and let }$\mathbb{L}:=%
\mathbb{L}_{\mathcal{M}}=\mathbb{L}_{\mathcal{N}}.$ \textit{Then }$\mathcal{M%
}$\textit{\ and }$\mathcal{N}$\textit{\ are isomorphic if the following
three conditions are satisfied}:\textit{\ }\medskip

\noindent $(I)$ $\mathrm{Cod}_{W}(\mathcal{M})=\mathrm{Cod}_{W}(\mathcal{N})$%
.\medskip

\noindent $(II)$ $\mathrm{Th}_{\mathbb{L}}(\mathcal{M)}=\mathrm{Th}_{\mathbb{%
L}}\mathrm{(}\mathcal{N)}$\textit{.}\medskip

\noindent $(III)$ \textit{For some \textbf{nonstandard} ordinals} $\tau _{M}$
\textit{of }$\mathcal{M}$ \textit{and }$\tau _{N}$ \textit{of }$\mathcal{N}$
\textit{there are} $S_{M}\subseteq M$ \textit{and }$S_{N}\subseteq N$
\textit{such that} $S_{M}$ \textit{is a separative }$\tau _{M}$\textit{%
-satisfaction class} \textit{over} $\mathcal{M}$, \textit{and} $S_{N}$
\textit{is a separative} $\tau _{N}$\textit{-satisfaction class} \textit{over%
} $\mathcal{N}$.\medskip

\noindent \textbf{Proof.}~The isomorphism between $\mathcal{M}$ and $%
\mathcal{N}$ can be built by a routine back-and-forth construction once we
establish Claim 3.1.1 below, for which we introduce the following convention:

\begin{itemize}
\item Given an $n$-tuple $\overline{a}=(a_{0},\cdot \cdot \cdot ,a_{n-1})$
from\textit{\ }$M$ (where $n\in \omega $), and an $n$-tuple $\overline{b}%
=(b_{0},\cdot \cdot \cdot ,b_{n-1})$ from\textit{\ }$N$, we write $\overline{%
a}\thicksim \overline{b}$ as a shorthand for the following statement, where $%
\overline{x}$ is an $n$-tuple:
\end{itemize}

\begin{center}
for each $n$-ary formulae $\varphi (\overline{x})$ of $\mathbb{L},$ $\varphi
(\overline{a})\in S_{M}$ iff $\varphi (\overline{b})\in S_{N}.$
\end{center}

\noindent Note that by Proposition 2.3, $\overline{a}\thicksim \overline{b}$
iff for all $n$-ary formulae $\varphi (\overline{x})$ of $\mathbb{L},$ $%
\mathcal{M}\models \varphi (\overline{a})$ iff $\mathcal{N}\models \varphi (%
\overline{a})$. In particular, $\varnothing \thicksim \varnothing ,$ by
condition $(II)$ of the lemma, where $\varnothing $ is the \textquotedblleft
0-tuple", i.e., the empty sequence.\medskip

\noindent \textbf{3.1.1.~Claim}.~\textit{Suppose} $\overline{a}\thicksim
\overline{b}.$ \textit{Then}:$\medskip $

\noindent $(i)$ \textit{For every} $a\in M$ \textit{there is} \textit{some} $%
b\in N$ \textit{such that} $\left( \overline{a},a\right) \thicksim \left(
\overline{b},b\right) .\medskip $

\noindent $(ii)$ \textit{For every} $b\in N$ \textit{there is} \textit{some}
$a\in M$ \textit{such that} $\left( \overline{a},a\right) \thicksim \left(
\overline{b},b\right) .$\medskip

\noindent By symmetry it suffices to verify part $(i)$ of Claim 3.1.1.
Observe that since $\mathcal{M}$ and $\mathcal{N}$ share the same
well-founded part $W$, we can fix an ordinal $\eta $ such that $\eta =%
\mathrm{o}(\mathcal{M})=\mathrm{o}(\mathcal{N}),$ and $N(\alpha )=M(\alpha )$
for all $\alpha <\eta.$

\medskip

\noindent Given $a\in M$, let

\begin{center}
$X:=\{\varphi (\overline{v},v):\varphi (\overline{v},v)$ is an $\left(
n+1\right) $-ary formula of $\mathbb{L}$, and $\varphi (\overline{a},a)\in
S_{M}\}.$
\end{center}

\noindent A routine argument shows that $X\in \mathrm{Cod}_{W}(\mathcal{M})$
(using the assumption that $S_{M}$ is a separative $\tau _{M}$-satisfaction
class over $\mathcal{M}$ and $\tau _{M}$ is a nonstandard ordinal of $%
\mathcal{M}$). So by assumption $(I)$ of Lemma 3.1, $X\in \mathrm{Cod}_{W}(%
\mathcal{N}).$ Hence there is some $c\in M$ such that $X=W\cap \mathrm{Ext}_{%
\mathcal{N}}(c)$. For any $\alpha \in \mathrm{Ord}^{\mathcal{M}}$, consider
the elements $c_{\alpha }$ and $d_{\alpha }$ of $N$, such that the following
holds in $\mathcal{N}$:
\begin{equation*}
c_{\alpha }=\{x\in c:x\in \mathrm{V}(\alpha )\}\ \mathrm{and}\ d_{\alpha
}=\{x\in \mathrm{V}(\alpha ):x\notin c\}.
\end{equation*}%
Then for each $\alpha <\eta $, both $c_{\alpha }$ and $d_{\alpha }\in W.$
Also, in light of our convention of identifying $W$ with its transitive
collapse, for each $w\in W$ we have:

\begin{center}
$w=$ $\mathrm{Ext}_{\mathcal{M}}(w)=\mathrm{Ext}_{\mathcal{N}}(w)$.
\end{center}

\noindent The choice of $c_{\alpha }$ and $d_{\alpha }$ together with the
compositional properties of $S_{M}$ allows us to conclude:\medskip

\noindent (1) For all $\alpha \in \eta \overset{\psi _{\alpha }(\overline{a})%
}{~\overbrace{\left( \exists x\left( \left( \bigwedge\limits_{\varphi (%
\overline{v},v)\in c_{\alpha }}\varphi (\overline{a},x)\right) \wedge \left(
\bigwedge\limits_{\varphi (\overline{v},v)\in d_{\alpha }}\lnot \varphi (%
\overline{a},x)\right) \right) \right) }}\in S_{M} $.\medskip

\noindent Observe that $\psi _{\alpha }(\overline{x})$ is a formula of $%
\mathbb{L}.$ Putting (1) together with the assumption $\overline{a}\thicksim
\overline{b}$ yields $\psi _{\alpha }(\overline{b})\in S_{N}$, i.e., \medskip

\noindent (2) For all $\alpha \in \eta ~\left( \exists x\left( \left(
\bigwedge\limits_{\varphi (\overline{v},v)\in c_{\alpha }}\varphi (\overline{%
b},x)\right) \wedge \left( \bigwedge\limits_{\varphi (\overline{v},v)\in
d_{\alpha }}\lnot \varphi (\overline{b},x)\right) \right) \right) \in S_{N}$%
.\medskip

\noindent The key observation at this point is that there is a first order
formula $\theta (\mathrm{S},x,y,\overline{z})$ in the language of set theory
augmented with the predicate \textrm{S} such that (2) can be re-expressed
as:\medskip

\noindent (3) For all $\alpha \in \eta ,\ (\mathcal{N},S_{N})\models \theta (%
\mathrm{S},\alpha ,c,\overline{b}).$\medskip

\noindent By invoking Overspill (Proposition 2.6) in the expanded structure $%
(\mathcal{N},S_{N})$, there is some \textit{nonstandard} ordinal $\gamma $
of $\mathcal{N}$ such $(\mathcal{N},S_{N})\models \theta (\mathrm{S,}\gamma
,c,\overline{b})$, i.e., \medskip

\noindent (4) $(\mathcal{N},S_{N})\models \mathrm{S}\left( \exists x\left(
\left( \bigwedge\limits_{\varphi (\overline{v},v)\in c_{\gamma }}\varphi (%
\overline{b},x)\right) \wedge \left( \bigwedge\limits_{\varphi (\overline{v}%
,v)\in d_{\gamma }}\lnot \varphi (\overline{b},x)\right) \right) \right) .$%
\medskip

\noindent By coupling (4) together with the assumption that $\left( \mathcal{%
N},S_{N}\right) $ satisfies condition $(IV)$ of $\mathrm{Sat}(\tau_{N} ,S)$
(as in Definition 2.2), the existential statement deemed true in (4) by the
interpretation $S_{N}$ of $\mathrm{S}$ is witnessed by some $b\in N.$ It
should be clear that this is the desired element $b\in N,$ i.e., $\left(
\overline{a},a\right) \thicksim \left( \overline{b},b\right) .$ This
concludes the proof of Claim 3.1.1, and therefore of Lemma 3.1.\hfill $%
\square $\medskip

We now present an easy lemma (Lemma 3.2), and an old theorem of Hutchinson
(Theorem 3.3); they will allow us to arrange the hypotheses of Lemma 3.1 in
the proof of Theorem A. \medskip

\noindent \textbf{3.2.~Lemma.}~(ZFC) \textit{Let} $\lambda $ \textit{be a
strongly inaccessible cardinal, }$S\subseteq \mathrm{V}(\lambda )$, \textit{%
and let }

\begin{center}
$C:=\left\{ \delta <\lambda :(\mathrm{V}(\delta ),\in ,S\cap \mathrm{V}%
(\delta ))\prec (\mathrm{V}(\lambda ),\in ,S)\right\} $.
\end{center}

\noindent \textit{Then} $C$ \textit{is closed and unbounded in} $\lambda .$%
\medskip

\noindent \textbf{Proof.}~$C$ is clearly closed by the elementary chain
theorem, so we will concentrate on demonstrating the unboundedness of $C$.
Fix a well-ordering $\vartriangleleft $ of $\mathrm{V}(\lambda )$, and for
any $A\subseteq \mathrm{V}(\lambda )$, let $\mathcal{H}(A)$ be the submodel
of $\left( \mathrm{V}(\lambda ),\in ,S\right) $ whose universe $H(A)$
consists of the elements of $\mathrm{V}(\lambda )$ that are first order
definable in the expanded structure $\left( \mathrm{V}(\lambda ),\in
,S,\vartriangleleft ,a\right) _{a\in A}.$ Clearly $\left\vert
H(A)\right\vert =\min \{\aleph _{0},\left\vert A\right\vert \}$, and by
Tarski's test $\mathcal{H}(A)\prec $ $\left( \mathrm{V}(\lambda ),\in
,S\right) .$ Given an ordinal $\alpha <\lambda ,$ we will exhibit $\beta $
such that $\alpha \leq \beta \in C.$ To this end, consider the sequence of
models $\left\langle \mathcal{M}_{n}:n\in \omega \right\rangle $ and
sequence of ordinals $\left\langle \mathcal{\alpha }_{n}:n\in \omega
\right\rangle $ defined by the following recursive clauses:

\begin{itemize}
\item $\mathcal{\alpha }_{0}:=\alpha ,$ and $\mathcal{M}_{0}:=\mathcal{H}%
(\alpha _{0}).$

\item $\alpha _{n+1}:=\sup \{\beta <\kappa :\beta \in M_{n}\}$, and $%
\mathcal{M}_{n+1}:=\mathcal{H}(\mathrm{V}(\alpha _{n+1})).$
\end{itemize}

\noindent The strong inaccessibility of $\lambda $ guarantees that $\mathcal{%
M}_{n}$ and $\alpha _{n}$ are well-defined for each $n\in \omega $, and that
$\left\{ \mathcal{\alpha }_{n}:n\in \omega \right\} $ is bounded in $\lambda
$. Let $\beta :=\sup \left\{ \mathcal{\alpha }_{n}:n\in \omega \right\} .$
It is routine to verify that $\alpha \leq \beta \in C.$\hfill $\square $%
\medskip

The following theorem was established by Hutchinson \cite{Hutchinson} using
the omitting types theorem. As shown in \cite[Theorem 2.12]{Ali-TAMS} this
result can also be proved for models of ZFC using generic
ultrapowers.\medskip

\noindent \textbf{3.3.~Theorem.}~(Hutchinson)\textbf{\ }\textit{Suppose }$%
\lambda $ \textit{is a regular cardinal in a countable model} $\mathcal{K}$
\textit{of} ZF. \textit{Then there is an elementary extension} $\mathcal{K}%
^{\ast }$ \textit{of} $\mathcal{K}$ \textit{satisfying the following two
properties}:\medskip

\noindent \textbf{(a)}\textit{\ }$\mathcal{K}^{\ast }$ \textit{does not}
\textit{\textquotedblleft perturb\textquotedblright\ any ordinal of} $%
\mathcal{K}$ \textit{that is below} $\lambda ,$ \textit{i.e., if }$\mathcal{K%
}\models \alpha \in \lambda $, \textit{then} $\mathrm{Ext}_{\mathcal{K}%
}(\alpha )=\mathrm{Ext}_{\mathcal{K}^{\ast }}(\alpha )$.\textit{\ \smallskip
}

\noindent \textbf{(b) }$\mathrm{Ext}_{\mathcal{K}^{\ast }}(\lambda )\
\backslash \ \mathrm{Ext}_{\mathcal{K}}(\lambda )$, \textit{when ordered} by
$\in ^{\mathcal{K}^{\ast }}$, \textit{has no first element }(\textit{under
the ordering }$\in ^{\mathcal{K}^{\ast }}$)\textit{.}\medskip

\noindent \textbf{3.4.~Remark.}~Condition (a) of Theorem 3.3 ensures that if
$k\in K$ and $\mathcal{K}\models \left\vert k\right\vert <\lambda $, then $%
\mathcal{K}$ does not perturb $k$. To see this, choose $f$ and $\alpha $ in $%
K$ such that:

\begin{center}
$\mathcal{K}\models $ \textquotedblleft $\alpha \in \lambda $ and $f:\alpha
\rightarrow k$ and $f$ is a bijection\textquotedblright .
\end{center}

\noindent Then since $\mathcal{K}\prec \mathcal{K}^{\ast },$ $\mathrm{Ext}_{%
\mathcal{K}^{\ast }}(k)=\{f(x)\in K^{\ast }:x\in \mathrm{Ext}_{\mathcal{K}%
^{\ast }}(\alpha )\},$ and so together with the assumption $\mathrm{Ext}_{%
\mathcal{K}^{\ast }}(\alpha )=\mathrm{Ext}_{\mathcal{K}}(\alpha )$, this
makes it clear that $\mathrm{Ext}_{\mathcal{K}^{\ast }}(k)=\mathrm{Ext}_{%
\mathcal{K}}(k).$ Therefore, if $\mathcal{K}$ is well-founded, and $\lambda $
is strongly inaccessible in $\mathcal{K}$, $\mathrm{WF}(\mathcal{K}^{\ast
})=K(\lambda )=\mathrm{WF}(\mathcal{K}^{\ast }(\lambda )).$\medskip

We are now ready to present the proof of Theorem A. Recall that $\mathcal{M}$
is \textit{definably well-founded} if every element of $\mathcal{M}$\ that
is first order definable in $\mathcal{M}$\textit{\ }(without parameters)%
\textit{\ }is in the well-founded part of\textit{\ }$\mathcal{M}$. In
particular, if $\mathcal{M}$ is definably well-founded, then for any
parameter-free definable $\alpha \in \mathrm{Ord}^{\mathcal{M}},$ the
predecessors of $\alpha $ form a well-ordered set as viewed externally, and
thus a nonstandard ordinal of $\mathcal{M}$ (if any) dwarfs any definable
ordinal of $\mathcal{M}$. In particular a definably well-founded model is an
$\omega $-model. We say that $\mathcal{M}$ is \textit{cofinally condensable}
if the collection of $\alpha \in \mathrm{Ord}^{\mathcal{M}}$ such that $%
\mathcal{M}\cong \mathcal{M(\alpha )}\prec _{\mathbb{L}_{\mathcal{M}}}%
\mathcal{M}$ is unbounded in\textit{\ }$\mathrm{Ord}^{\mathcal{M}}$.\medskip

\noindent \textbf{Theorem A.}~\textit{Assuming that }$\mathit{\mathrm{ZFC}}$%
\textit{\ }+\textit{\ \textquotedblleft there exists an inaccessible
cardinal\textquotedblright\ has a well-founded model, there is a model }$%
\mathcal{M}$ \textit{of} $\mathrm{ZFC}$\textit{\ that is both definably
well-founded and cofinally condensable.\ }\medskip

\noindent \textbf{Proof.}~The proof is carried out in two steps, the first
takes place within an appropriately chosen model $\mathcal{K}$ of $\mathrm{%
ZFC}$; the second step is performed outside of $\mathcal{K}$.\medskip

\noindent \textbf{Step 1.}~If the theory $\mathrm{ZFC}$ + \textquotedblleft
there exists an inaccessible cardinal" has a well-founded model, then by the
L\"{o}wenheim-Skolem theorem and the fact that ZF proves that GCH holds in
the constructible universe, there is a countable well-founded model that
contains a strongly inaccessible cardinal (since if $\lambda $ is
inaccessible in a model $\mathcal{K}$, then $\lambda $ is also inaccessible
in the constructible universe $\mathrm{L}^{\mathcal{K}}$ of $\mathcal{K}$;
and under GCH every inaccessible cardinal is strongly inaccessible). Let $%
\mathcal{K}$ be a countable well-founded model that contains a
\textquotedblleft cardinal\textquotedblright\ $\lambda $ that is strongly
inaccessible in the sense of $\mathcal{K}$. By collapsing $\mathcal{K}$ we
may assume that $\mathcal{K}=(K,\in ).$ By Remark 2.4 and Lemma 3.2 we can
get hold of elements $s$ and $u$ of $K$ satisfying the following conditions:$%
\medskip $

\noindent $(i)$ $\mathcal{K}\models $ \textquotedblleft $s$\textit{\ }is a
separative $\infty $-satisfaction class for $(\mathrm{V}(\lambda ),\in )$"$%
.\medskip $

\noindent $(ii)$ $\mathcal{K}\models $ \textquotedblleft $u$ is unbounded in
$\lambda $ and $\forall \delta \in u$ $(\mathrm{V}(\delta ),\in ,s\cap
\mathrm{V}(\delta ))\prec (\mathrm{V}(\lambda ),\in ,s)$\textquotedblright .$%
\medskip $

\noindent \textbf{Step 2.}~By Theorem 3.3 and Remark 3.4, there is an
elementary extension $\mathcal{K}^{\ast }$ of $\mathcal{K}$ such that $%
\mathrm{WF}(\mathcal{K}^{\ast })=K(\lambda )=\mathrm{WF}(\mathcal{K}^{\ast
}(\lambda )).$ \medskip

We claim that $\mathcal{K}^{\ast }(\lambda )$ is definably well-founded and
cofinally condensable. $\mathcal{K}^{\ast }(\lambda )$ is definably
well-founded since $\mathcal{K}(\lambda )\prec\mathcal{K}^{\ast }(\lambda )$%
, and $\mathcal{K}(\lambda )$ is well-founded. Recall that $\mathcal{K}$
thinks that $u$ is an unbounded subset of $\lambda $. Since $\mathcal{K}%
\prec \mathcal{K}^{\ast }$, to verify that $\mathcal{K}^{\ast }(\lambda )$
is cofinally condensable it suffices to show that if $\mathcal{K}^{\ast
}\models \delta \in u$, and $\delta \in K^{\ast }\backslash K$ \
(equivalently: $\delta $ is a nonstandard element of $\mathrm{Ext}_{\mathcal{%
K}^{\ast }}(u ))$, then $\mathcal{K}^{\ast }(\lambda )\cong \mathcal{K}%
^{\ast }(\delta ).$ This is precisely where Lemma 3.1 comes into the
picture. If $S:=\mathrm{Ext}_{\mathcal{K}^{\ast }}(s)$, then $(i)$ and $(ii)$
assure us that the assumptions of Lemma 3.1 are satisfied if we choose any
nonstandard ordinal $\tau_{0}$ below $\delta$ and let:

\begin{center}
$\mathcal{M}:=\mathcal{K}^{\ast }(\lambda ),$ $\mathcal{N}:=\mathcal{K}%
^{\ast }(\delta );$ $\tau _{M}:=\tau_{0} $, $\tau _{N}:=\tau_{0}$, $S_{M}:=S$
and $S_{N}:= S\cap K^{\ast }(\delta )$.
\end{center}

\noindent Hence by Lemma 3.1, $\mathcal{K}^{\ast }(\lambda )\cong \mathcal{K}%
^{\ast }(\delta )$, thus concluding the proof of Theorem A. \hfill $\square $%
\medskip

\noindent \textbf{3.5. Remark.}~After seeing Theorem A Corey Switzer asked
the author whether there are\textit{\ uncountable} definably well-founded
cofinally condensable models of set theory. Assuming the existence of a
weakly compact cardinal, the answer is in the positive. Here we outline the
construction of such a model. Let $\lambda $ be a weakly compact cardinal.
Then if $T$ is a theory\ formulated in the infinitary logic $\mathbb{L}%
_{\lambda ,\lambda }$ such that $\left\vert T\right\vert =\lambda $, and
every subset of $T$ of cardinality less than $\lambda $ has a model, then $T$%
\ has a model. This property of $\lambda $ can be used to show that if $%
X\subseteq \mathrm{V}(\lambda ),$ then the model $(\mathrm{V}(\lambda ),\in
,X)$ has an $\mathbb{L}_{\lambda ,\lambda }$-elementary end extension $(%
\mathcal{K},X_{K})$ such that $\mathrm{Ord}^{\mathcal{K}}\backslash \lambda $
has no least ordinal, and thus $\mathrm{WF}(\mathcal{K)}=\mathrm{V}(\lambda
) $. In particular, $\mathcal{K}$ is definably well-founded. In order to
ensure that $\mathcal{K}$\ is also cofinally condensable, we can choose $X$
to be the $\mathbb{L}_{\lambda ,\lambda }$ satisfaction class for $\left(
\mathrm{V}(\lambda ),\in \right) $. We then prove a suitable adaptation of
Lemma 3.1, where the countability assumption of Lemma 3.1 is replaced with
the assumption that both $\mathcal{M}$ and $\mathcal{N}$ have cardinality $%
\lambda $ for some strongly inaccessible $\lambda $, $\mathrm{WF}(\mathcal{M}%
)=\mathrm{WF}(\mathcal{N})=\mathrm{V}(\lambda )$, and there are $%
S_{M}\subseteq M$, $\tau _{M}\in \mathrm{Ord}^{\mathcal{M}},$ $%
S_{N}\subseteq M$, and $\tau _{N}\in \mathrm{Ord}^{\mathcal{N}}$ such that $%
\tau _{M}$ and $\tau _{N}$ are respectively nonstandard\textit{\ cardinals}
of $\mathcal{M}$ and $\mathcal{N}$, $S_{M}$ satisfies Tarski's compositional
axioms over $\mathcal{M}$ for all formulae of the logic $\mathbb{L}_{\delta
,\delta }^{\mathcal{M}}$ for $\delta =\tau _{M}$, and $S_{N}$ is satisfies
Tarski's compositional axioms over $\mathcal{N}$ for all formula of the
logic $\mathbb{L}_{\delta ,\delta }^{\mathcal{N}}$ for $\delta =\tau _{N}$.
Thus $S_{M}$ correctly calculates the truth value of all $\mathbb{L}%
_{\lambda ,\lambda }$ sentences over $\mathcal{M}$, and $S_{N}$ correctly
calculates the truth value of all $\mathbb{L}_{\lambda ,\lambda }$ sentences
over $\mathcal{N}$. Using an argument very similar to the proof of Theorem
A, one then shows that $\mathcal{K}$\ is also cofinally condensable. \bigskip

\begin{center}
\textbf{4.~PROOF OF THEOREM B}\bigskip
\end{center}

We first lay out a series of definitions and lemmas before presenting the
proof of Theorem B.

\begin{itemize}
\item Throughout the section, $\mathcal{M}$ is assumed to be a nonstandard
model of $\mathrm{ZF}$ and $W:=\mathrm{WF}(\mathcal{M)}$.$\medskip $
\end{itemize}

\noindent \textbf{4.1.~Definition.}~A structure $\mathcal{N}$ is $\mathrm{Cod%
}_{W}(\mathcal{M})$-\textit{saturated} if for every type $p(x,y_{1},\cdot
\cdot \cdot ,y_{k})$, and for every $k$-tuple $\overline{a}$ of parameters
from $\mathcal{N}$, $p(x,\overline{a})$ is realized in $\mathcal{N}$
provided the following three conditions are satisfied:$\medskip $

\noindent $(i)$ $p(x,\overline{y})\subseteq \mathbb{L}_{\mathcal{M}}$.$%
\medskip $

\noindent $(ii)$ $p(x,\overline{y})\in \mathrm{Cod}_{W}(\mathcal{M})$.$%
\medskip $

\noindent $(iii)$ $\forall w\in W\ \mathcal{N}\models \exists x\left(
\bigwedge\limits_{\varphi \in p(x,\overline{y})\cap w}\varphi (x,\overline{a}%
)\right) .$

\begin{itemize}
\item In the interest of concision, we say that $\mathcal{M}$\ is $W$-%
\textit{saturated} if $\mathcal{M}$ is $\mathrm{Cod}_{W}(\mathcal{M})$%
-saturated.
\end{itemize}

\noindent \textbf{Remark 4.1.1.}~It is not hard to see that if $\mathcal{M}$
is $\omega $-nonstandard, then $\mathcal{M}$ is $W$-saturated iff $\mathcal{M%
}$ is recursively saturated. We should also point out that a notion closely
related to $\mathrm{Cod}_{W}(\mathcal{M})$-saturation was introduced in
Ressayre's paper \cite{Ressayre-APAL} (dubbed $\alpha $-recursive
saturation) where it was used as a tool for studying the model theory of
admissible languages, as well as certain aspects of descriptive set theory
(see also Barwise \cite[p.143]{Barwise-book}, Schlipf \cite[p.164-165]%
{Schlipf-APAL}). A trick similar to the one that shows that recursive
saturation coincides with $W$-saturation for $\omega $-nonstandard models
can be used to show that, more generally, if $\mathcal{M}$ is nonstandard,
then $\mathcal{M}$\ is $W$-\textit{saturated} iff $\mathcal{M}$ is $\mathrm{o%
}(\mathcal{M})$-recursively saturated. \medskip

\noindent \textbf{4.2.~Lemma.}~\textit{If} $\gamma $ \textit{is a limit
ordinal of} $\mathcal{M}$ \textit{and} $\gamma $ \textit{is nonstandard},
\textit{then} $\mathcal{M}(\gamma )$ \textit{is} $W$-\textit{saturated.}%
\medskip

\noindent \textbf{Proof.}~Given a 1-type $p(x,\overline{a})$, where $%
\overline{a}$ is a $k$-tuple of parameters from $\mathcal{M}(\gamma )$ such
that conditions $(i)$, $(ii)$, and $(iii)$ of Definition 4.1 hold, choose $%
c\in M(\gamma )$ such that $p(x,\overline{y})=\mathrm{Ext}_{\mathcal{M}%
(\gamma )}(c)\cap W,$ and let $\theta (z,\overline{y})$ be the following
formula in the language of set theory augmented with a predicate $\mathrm{S}$%
:

\begin{center}
$\exists x\left[ \forall \varphi (x,\overline{y})\in c\cap \mathrm{V}(z)\
\mathrm{S}(\varphi (c_{x},c_{y_{1}},\cdot \cdot \cdot ,c_{y_{k}})\right] .$
\end{center}

\noindent By Remark 2.4 there is some $s\in M$ such that $S:=\mathrm{Ext}_{%
\mathcal{M}}(s)$ is a separative $\gamma $-satisfaction class on $\mathcal{M}%
(\gamma ).$ Since for all $\alpha \in \mathrm{o}(\mathcal{M})$, $\left(
\mathcal{M}(\gamma ),S\right) $ satisfies $\theta (\alpha ,\overline{a})$ by
Proposition 2.3, $\theta (\gamma ^{\prime },\overline{a})$ holds in $\left(
\mathcal{M}(\gamma ),S\right) $ for some $\gamma ^{\prime }\in \mathrm{Ord}^{%
\mathcal{M}(\gamma )}\backslash W$ by Proposition 2.6 (Overspill), which
makes it evident that $p(x,\overline{a})$ is realized in $\mathcal{M}(\gamma
).$ Note that a slight modification of the proof shows that, more generally,
any structure that \textquotedblleft lives\textquotedblright\ in $\mathcal{M}
$ (i.e., has an isomorphic copy that is coded in $\mathcal{M}$) is $\mathrm{%
Cod}_{W}(\mathcal{M})$-saturated.\hfill $\square $\medskip

\noindent \textbf{4.3.~Lemma.}~\textit{Given countable nonstandard models} $%
\mathcal{M}$ \textit{and} $\mathcal{N}$ \textit{of} $\mathrm{ZF}$, $\mathcal{%
M}\cong \mathcal{N}$ \textit{provided the following two conditions hold}%
:\medskip

\noindent \textbf{(a)} $\mathcal{M}$\textit{\ and }$\mathcal{N}$\textit{\
have the same well-founded part }$W$\textit{, }$\mathrm{Cod}_{W}(\mathcal{M}%
)=\mathrm{Cod}_{W}(\mathcal{M}),$ and $\mathrm{ZF}(\mathbb{L})\subseteq
\mathrm{Th}_{\mathbb{L}}\mathrm{(}\mathcal{M)}=\mathrm{Th}_{\mathbb{L}}%
\mathrm{(}\mathcal{N)}$ \textit{for} $\mathbb{L}:=\mathbb{L}_{\mathcal{M}}=%
\mathbb{L}_{\mathcal{N}}$. \medskip

\noindent \textbf{(b)} \textit{Both} $\mathcal{M}$ \textit{and} $\mathcal{N}$
\textit{are} $W$-\textit{saturated.}\medskip

\noindent \textbf{Proof.} ~This lemma is a distillation of Lemma 3.1 since
the proof of Claim 3.1.1 can be readily modified to show that Claim 3.1.1
holds with the assumptions of Lemma 4.3 once we make the observation that
the $W$-saturation of $\mathcal{M}$ implies that the $\mathbb{L}$-type of
any finite tuple in $\mathcal{M}$ is a member of $\mathrm{Cod}_{W}(\mathcal{M%
})$ (and of course the same goes for $\mathcal{N}$). To verify this
observation, first consider the following type $p(x,\overline{y}):$

\begin{center}
$p(x,\overline{y}):=\left\{ \varphi (\overline{y})\leftrightarrow \left(
\varphi (\overline{y})\in x\right) :\varphi (y_{1},\cdot \cdot \cdot
,y_{n})\in \mathbb{L}\right\} .$
\end{center}

\noindent It is easy to see that $p(x,\overline{y})\in \mathrm{Cod}_{W}(%
\mathcal{M})$. Given $\overline{a}\in M^{n}$ and $\alpha \in \mathrm{o}(%
\mathcal{M})$, for sufficiently large $\beta \in \mathrm{o}(\mathcal{M}),$
we have:

\begin{center}
For all $\varphi (\overline{y})\in M(\alpha )$, $\mathcal{M}\models \varphi (%
\overline{a})$ iff $\mathcal{M}\models \mathrm{Sat}_{\beta }(\varphi (%
\overline{a})).$
\end{center}

\noindent Together with Proposition 2.5 and the assumption that $\mathcal{M}%
\models \mathrm{ZF}(\mathbb{L})$ we conclude that for each $\alpha \in
\mathrm{o}(\mathcal{M})$ the set $\left\{ \varphi (\overline{y})\in M(\alpha
):\mathcal{M}\models \varphi (\overline{a})\right\} $ is coded in $\mathcal{M%
}$. This makes it evident that the three conditions of Definition 4.1 are
met and therefore by the assumption of $W$-saturation of $\mathcal{M}$,
there is an element $c\in M$ such that for all $n$-ary $\mathbb{L}$-formulae
$\varphi (\overline{y})\in \mathbb{L}$, we have:

\begin{center}
$\mathcal{M}\models \varphi (\overline{a})$ iff $\varphi (\overline{y})\in
\mathrm{Ext}_{\mathcal{M}}(c),$
\end{center}

\noindent which shows that the $\mathbb{L}$-type of $\overline{a}$ is a
member of $\mathrm{Cod}_{W}(\mathcal{M}).$\hfill $\square $\medskip

\noindent \textbf{4.4. Lemma.}~\textit{Suppose there is an unbounded
collection of} $\alpha \in \mathrm{Ord}^{\mathcal{M}}$ \textit{such that} $%
\mathcal{M(\alpha )}\prec _{\mathbb{L}_{\mathcal{M}}}\mathcal{M}$. \textit{%
Then} $\mathcal{M}$ \textit{is} $W$-\textit{saturated.}\medskip

\noindent \textbf{Proof.}~This directly follows from Lemma 4.2 and
Proposition 2.7 (Elementary Chains). $\hfill \square $\medskip

\noindent \textbf{4.5. Lemma. }\textit{If there is some }$\alpha \in \mathrm{%
Ord}^{\mathcal{M}}$ \textit{with }$\mathcal{M(\alpha )}\prec _{\mathbb{L}_{%
\mathcal{M}}}\mathcal{M}$, \textit{then} $\mathcal{M}\models \mathrm{ZF}(%
\mathbb{L}_{\mathcal{M}})$.\medskip

\noindent \textbf{Proof.} Since $\mathcal{M}$ is assumed to be a model of $%
\mathrm{ZF}$, we just need to verify that $\mathcal{M}$ satisfies $\mathrm{%
Sep}(\mathbb{L}_{\mathcal{M}})$ and $\mathrm{Coll}(\mathbb{L}_{\mathcal{M}%
}). $ In light of the assumption that $\mathcal{M(\alpha )}$ is an $\mathbb{L%
}_{\mathcal{M}}$-elementary submodel of $\mathcal{M}$, it suffices to verify
that $\mathcal{M}(\alpha )$ satisfies $\mathrm{Sep}(\mathbb{L}_{\mathcal{M}%
}) $ and $\mathrm{Coll}(\mathbb{L}_{\mathcal{M}})$. To see that $\mathcal{%
M(\alpha )}\models \mathrm{Sep}(\mathbb{L}_{\mathcal{M}})$, suppose $\psi
(x)\in \mathbb{L}_{\mathcal{M}}$ (where $\psi (x)$ is allowed to have
parameters from $M(\alpha ))$, and $m\in M(\alpha ).$ Consider

\begin{center}
$K:=\left\{ k\in M:\mathcal{M(\alpha )}\models k\in m\wedge \psi (k)\right\}
.$
\end{center}

\noindent By Remark 2.4, there is a separative $\alpha $-satisfaction class $%
S$ on $\mathcal{M(\alpha )}$. Let

\begin{center}
$K^{\prime }:=\left\{ k\in M:\left( \mathcal{M(\alpha )},S\right) \models
k\in m\wedge \mathrm{S}(\psi (c_{k}))\right\} $.
\end{center}

\noindent By Proposition 2.3, $K=K^{\prime }$. On the other hand, since $S$
is separative, $K^{\prime }$ is coded in $\mathcal{M(\alpha )}$. This
concludes the proof that $\mathrm{Sep}(\mathbb{L}_{\mathcal{M}})$ holds in $%
\mathcal{M}(\alpha )$. To verify that $\mathcal{M}(\alpha )$ satisfies $%
\mathrm{Coll}(\mathbb{L}_{\mathcal{M}}),$ suppose for some $\varphi (x,y)\in
\mathbb{L}_{\mathcal{M}}$ and for some $m$ in $M\mathcal{(\alpha )}$ we
have:\medskip

\noindent (1) $\mathcal{M(\alpha )}\models \forall x\in m$ $\exists y$ $%
\varphi (x,y).$\medskip

\noindent We need to verify:\medskip

\noindent (2) $\mathcal{M(\alpha )}\models \exists z\ \forall x\in m$ $%
\exists y\in z$ $\varphi (x,y).$\medskip

\noindent Define $f(x)$ in $\mathcal{M}$ to be the unique ordinal $\alpha $
that satisfies $\psi (x,\alpha ),$ where

\begin{center}
$\psi (x,\alpha ):=$ $\exists y\in \mathrm{V}(\alpha )$ $\left( \varphi
(x,y)\wedge \left( \forall \beta \in \alpha \ \forall y\in \mathrm{V}(\alpha
)\ \lnot \varphi (x,y)\right) \right) .$
\end{center}

\noindent Note that by (1) $f$ is well-defined in $\mathcal{M}$ for all $%
x\in m$. To establish (2) it suffices to show:\medskip

\noindent (3) $\mathcal{M(\alpha )}\models \exists \beta \in \mathrm{Ord\ }%
\forall x\in m$ $f(x)<\beta .$\medskip

\noindent Suppose (3) is false, then:\medskip

\noindent (4) $\mathcal{M(\alpha )}\models \forall \beta \in \mathrm{Ord\ }%
\exists x\in m$ $f(x)\geq \beta .$\medskip

\noindent So by the assumption $\mathcal{M(\alpha )}\prec _{\mathbb{L}_{%
\mathcal{M}}}\mathcal{M}$, (4) yields:\medskip

\noindent (5) $\mathcal{M}\models \forall \beta \in \mathrm{Ord\ }\exists
x\in m$ $f(x)\geq \beta .$\medskip

\noindent Pick $\beta \in \mathrm{Ord}^{\mathcal{M}}\backslash \ \mathrm{Ord}%
^{\mathcal{M(\alpha )}},$ then by (5) there is some $m_{0}\in \mathrm{Ext}_{%
\mathcal{M}}(m)$ ($=$ $\mathrm{Ext}_{\mathcal{M(\alpha )}}(m))$ such that $%
f(m_{0})\geq \beta $ holds in $\mathcal{M}$. This contradicts $\mathcal{%
M(\alpha )}\prec _{\mathbb{L}_{\mathcal{M}}}\mathcal{M}$ since $f$ is an $%
\mathbb{L}_{\mathcal{M}}$-definable function in $\mathcal{M}$, thereby
showing the veracity of (3). This concludes the verification of $\mathrm{Coll%
}(\mathbb{L}_{\mathcal{M}})$ in $\mathcal{M}.$\hfill $\square $\medskip

\noindent \textbf{4.6. Lemma.}~\textit{If }$\mathcal{M}$ \textit{is a
nonstandard model of }$\mathrm{ZF}(\mathbb{L}_{\mathcal{M}})$, \textit{and} $%
\mathcal{M}$ is $W$-\textit{saturated, then }$\mathcal{M}$\textit{\ is
cofinally condensable.}\medskip

\noindent \textbf{Proof.}~Fix any nonstandard $\gamma \in \mathrm{Ord}^{%
\mathcal{M}}$ and consider the type $p(x,\gamma )$ (where $\gamma $ is
treated as a parameter) consisting of the formula $\left( \gamma \in
x\right) \wedge \mathrm{Ord}(x)$ together with formulae of the form $\mathrm{%
Ref}_{\varphi }(x)$ (as in Proposition 2.8) as $\varphi $ ranges in $\mathbb{%
L}_{\mathcal{M}}$. It is easy to see that $p(x,y)$ satisfies conditions $(i)$
and $(ii)$ of Definition 4.1. Moreover, by Proposition 2.8 (Reflection) $%
p(x,\gamma )$ also satisfies condition $(iii)$ of Definition 4.1. Therefore
by the assumption of $W$-saturation of $\mathcal{M}$, $p(x,\mathcal{\gamma }%
) $ is realized in $\mathcal{M}$ by some $\gamma ^{\prime }$, which makes it
clear that $\gamma ^{\prime }$ is above $\gamma ,$ $\gamma ^{\prime }$ is a
nonstandard limit ordinal of $\mathcal{M}$, and $\mathcal{M(\gamma }^{\prime
}\mathcal{)}\prec _{\mathbb{L}_{\mathcal{M}}}\mathcal{M}$. Thanks to Lemmas
4.2 and 4.3, $\mathcal{M(\gamma }^{\prime }\mathcal{)}\cong \mathcal{M}$,
thus $\mathcal{M}$ is cofinally condensable.\hfill $\square $\medskip

We are now ready to establish Theorem B. The special case of Theorem B for
countable $\omega $-nonstandard models of $\mathrm{ZF}$ follows from
Schlipf's work on recursively saturated models of $\mathrm{ZF}$ in \cite%
{Schlipf-JSL} and \cite{Schlipf-PAMS}.\medskip

\noindent \textbf{Theorem B.}~\textit{The following are equivalent for a
countable model }$\mathcal{M}$ \textit{of }$\mathrm{ZF}$\textit{.}\medskip

\noindent \textbf{(a)} $\mathcal{M}$ \textit{is condensable}.\medskip

\noindent \textbf{(b) }$\mathcal{M}$ \textit{is cofinally condensable}.%
\textit{\ }\medskip

\noindent \textbf{(c)} $\mathcal{M}$ \textit{is nonstandard, and }$\mathcal{%
M(\alpha )}\prec _{\mathbb{L}_{\mathcal{M}}}\mathcal{M}$ \textit{for an
unbounded collection of }$\alpha \in \mathrm{Ord}^{\mathcal{M}}$.\medskip

\noindent \textbf{(d)} $\mathcal{M}$ \textit{is nonstandard and }$W$-\textit{%
saturated, and} $\mathcal{M}\models \mathrm{ZF}(\mathbb{L}_{\mathcal{M}}).$%
\medskip

\noindent \textbf{(e)} \textit{For some nonstandard ordinal} $\gamma $
\textit{of} $\mathcal{M}$ \textit{and some} $S\subseteq M$, $S$ \textit{is
an amenable}\footnote{%
This notion was defined in Definition 2.2.}\textit{\ }$\gamma $-\textit{%
satisfaction class on }$\mathcal{M}.$\medskip

\noindent \textbf{Proof.}~We will first show the equivalence of $(a)$, $(b)$%
, and $(c)$ by establishing $(b)\Rightarrow (a)\Rightarrow (c)\Rightarrow
(b) $. This will allow us to show the equivalence of $(d)$ with each of $(a)$%
, $(b)$, and $(c)$ by proving $(a)\Rightarrow (d)\Rightarrow (a)$. Finally,
we demonstrate $(a)\Rightarrow (e)\Rightarrow (d)$ to complete the proof.
\medskip

\noindent $\mathbf{(b)\Rightarrow (a)}.$ Trivial.\medskip

\noindent $\mathbf{(a)\Rightarrow (c)}.$ Suppose $\mathcal{M}$ is
condensable with $\mathcal{M}\cong \mathcal{M(\alpha )}\prec _{\mathbb{L}_{%
\mathcal{M}}}\mathcal{M}$. Then by \textquotedblleft
unwinding\textquotedblright\ the isomorphism between $\mathcal{M}$ and $%
\mathcal{M(\alpha )}$, we can readily obtain a sequence of models $%
\left\langle \mathcal{N}_{n}:n\in \omega \right\rangle $ such that $\mathcal{%
N}_{0}=\mathcal{M}$ and for all $n\in \omega $ the following hold:\medskip

\noindent (1) $\mathcal{N}_{n}=(\mathrm{V}(\alpha _{n}),\in )^{\mathcal{N}%
_{n+1}}$ for some $\alpha _{n}\in \mathrm{Ord}(\mathcal{N}_{n+1}).$\medskip

\noindent (2) $\mathcal{N}_{n}\prec _{\mathbb{L}}\mathcal{N}_{n+1}$, where $%
\mathbb{L}:=\mathbb{L}_{\mathcal{M}}$.\medskip

\noindent (3) $\mathcal{N}_{n}\cong \mathcal{M}$.\medskip

\noindent Let $\mathcal{N}:=\bigcup\limits_{n\in \omega }\mathcal{N}_{n}$.
By Proposition 2.7 (Elementary Chains), $\mathcal{N}_{n}\prec _{\mathbb{L}}%
\mathcal{N}$ for all $n\in \omega ,$ which together with Lemma 4.4 implies
that $\mathcal{N}$ is $W$-saturated, where $W=\mathrm{WF}(\mathcal{N})=%
\mathrm{WF}(\mathcal{N}_{n})$ for all $n\in \omega .$ By Lemma 4.2, $%
\mathcal{M}$ is also $W$-saturated. Therefore $\mathcal{N}\cong \mathcal{M}$
by Lemma 4.3, which in light of (2) and (3) and the unboundedness of $%
\{\alpha _{n}:n\in \omega \}$ in $\mathrm{Ord}^{\mathcal{N}}$ makes it clear
that (c) holds.\medskip

\noindent $\mathbf{(c)\Rightarrow (b)}.$ Assume (c). It is easy to see,
using Proposition 2.6 and Lemma 4.2, that $\mathcal{M}$ is $W$-saturated. By
(c) we can choose a nonstandard\textit{\ }$\gamma \in \mathrm{Ord}^{\mathcal{%
M}}$ arbitrarily high in \textrm{Ord}$^{\mathcal{M}}$ such that $\mathcal{M}%
_{\gamma }\prec _{\mathbb{L}_{\mathcal{M}}}\mathcal{M}$. Since by Lemma 4.2 $%
\mathcal{M}_{\gamma }$ is $W$-saturated, we can now invoke Lemma 4.3 to
conclude that $\mathcal{M\cong M}_{\gamma }\prec _{\mathbb{L}_{\mathcal{M}}}%
\mathcal{M}$, which makes it evident that (b) holds.\medskip

\noindent $\mathbf{(a)\Rightarrow (d)}.$ If $\mathcal{M\cong M(\alpha )}%
\prec _{\mathbb{L}_{\mathcal{M}}}\mathcal{M}$ for some $\alpha \in \mathrm{%
Ord}(\mathcal{M})$, then $\mathcal{M}$ is clearly nonstandard. Moreover, $%
\mathcal{M}$ is $W$-saturated by Lemma 4.2; and $\mathcal{M}$ satisfies $%
\mathrm{ZF}(\mathbb{L}_{\mathcal{M}})$ by Lemma 4.5.\medskip

\noindent $\mathbf{(d)\Rightarrow (a)}.$ This is justified by Lemma
4.6.\medskip

\noindent $\mathbf{(a)\Rightarrow (e)}.$ Suppose (a) holds and let $\alpha
\in \mathrm{Ord}(\mathcal{M})$ such that $\mathcal{M\cong M(\alpha )}\prec _{%
\mathbb{L}_{\mathcal{M}}}\mathcal{M}$. By Proposition 2.5 for each $\delta
\in \mathrm{o}(\mathcal{M}),$ there is some $S\in M$ such that $S$ is a $%
\delta $-satisfaction predicate over $\mathcal{M(\alpha )}$ that is
definable in $\mathcal{M(\alpha )}$ by an $\mathbb{L}_{\mathcal{M(\alpha )}}$%
-formula. Since we have verified that $(a)\Rightarrow (d)$, $\mathcal{%
M(\alpha )}$ satisfies $\mathrm{ZF}(\mathbb{L}_{\mathcal{M}})$, which
assures us that $S$ is an amenable $\delta $-satisfaction predicate over $%
\mathcal{M(\alpha )}$. Proposition 2.6 (Overspill) can be readily invoked
(applied to $\mathcal{M}$, rather than $(\mathcal{M},S)$) to show there is
some nonstandard $\gamma $ in $\mathcal{M}$ such that $\mathcal{M}$
satisfies \textquotedblleft there is an amenable $\gamma $-satisfaction
class over $(\mathrm{V}(\alpha ),\in )$\textquotedblright . In light of the
assumption that $\mathcal{M}\cong \mathcal{M(\alpha )}$, this shows that (e)
holds. \medskip

To carry out the overspill argument, we will distinguish between the case
when $\mathcal{M}$ is $\omega $-standard, and the case when $\mathcal{M}$ is
$\omega $-nonstandard. If $\mathcal{M}$ is $\omega $-standard, the overspill
argument succeeds smoothly since by routine absoluteness considerations, for
$s\in M$ and $S:=\mathrm{Ext}_{\mathcal{M}}(s)$, we have:

\begin{center}
$S$ is amenable over $\mathcal{M}(\alpha )$ iff $\mathcal{M}\models $
\textquotedblleft $s$ is amenable over $\left( \mathrm{V}(\alpha ),\in
\right) $\textquotedblright .
\end{center}

\noindent However, since the left-to-right direction of the above
equivalence can break down for $\omega $-nonstandard models (e.g., for
models of $\mathcal{M}$ of ZF that satisfy $\lnot \mathrm{Con(ZF)})$ we will
spell out the overspill argument for the case that $\mathcal{M}$ is $\omega $%
-nonstandard in more detail. It is worth pointing out that in this case (a)
implies that $\mathcal{M}$\ is recursively saturated, and by the
resplendence property of countable recursively saturated models one can
readily conclude that $\mathcal{M}$ carries an amenable $j$-satisfaction
class for some nonstandard $j\in \omega ^{\mathcal{M}}.$ However, the
overspill argument we present establishes (e) without the assumption of
countability of $\mathcal{M}$. Within $\mathcal{M}$, for each $i\in \omega ^{%
\mathcal{M}}$ let $\mathrm{Repl}_{i}\mathrm{(S)}$ consist of all instances
of the replacement scheme in the language \{$\in ,\mathrm{S}$\textrm{\}}
whose length is at most $i.$ Then define a subset $S$ of $\mathrm{V}(\alpha
) $ to be $i$-\textit{amenable} over $\left( \mathrm{V}(\alpha ),\in \right)
$ if $\left( \mathrm{V}(\alpha ),\in ,S\right) \models \mathrm{Repl}_{i}%
\mathrm{(S)}$. Since $\mathcal{M(\alpha )}$ is a model of ZF and for each
\textquotedblleft real world\textquotedblright\ natural number $n$, there is
an $n$-satisfaction class over $\mathcal{M(\alpha )}$ that is first order
definable in $\mathcal{M}$, we may conclude: \medskip

\noindent (1) For each $n\in \omega $ $\mathcal{M}\models $
\textquotedblleft $\exists s$ ($s$ is $n$-amenable over $\left( \mathrm{V}%
(\alpha ),\in \right) )$\textquotedblright .\medskip

\noindent Therefore by Overspill, there is some nonstandard $j\in \omega ^{%
\mathcal{M}}$ such that:\medskip

\noindent (2) $\mathcal{M}\models $ \textquotedblleft $\exists s$ ($s$ is $j$%
-amenable over $\left( \mathrm{V}(\alpha ),\in \right) )$\textquotedblright
.\medskip

\noindent Let $s\in M$ be a witness to the existential statement in (2) and $%
S:=\mathrm{Ext}_{\mathcal{M}}(s)$. It is evident that $S$ is an amenable $j$%
-satisfaction class over $\mathcal{M}(\alpha )$, as desired.\medskip

\noindent $\mathbf{(e)\Rightarrow (d)}.$ The $W$-saturation of $\mathcal{M}$
is readily verifiable with a reasoning very similar to the proof of Lemma
4.2. To see that $\mathcal{M}\models \mathrm{ZF}(\mathbb{L}_{\mathcal{M}}),$
it is sufficient to verify that the replacement scheme holds in $\mathcal{M}$%
\ for all $\mathbb{L}_{\mathcal{M}}$-formulae. To this end, let $\varphi
(x,y)$ be an $\mathbb{L}_{\mathcal{M}}$-formula (possibly with parameters
from $\mathcal{M}$), and suppose that $\mathcal{M}\models \forall x\exists
!y\varphi (x,y).$ Let $f:M\rightarrow M$ be the function whose graph is
described by $\varphi .$ Given $c\in M,$ we want to show that there is $d\in
M$ such that

\begin{center}
$\mathrm{Ext}_{\mathcal{M}}(d)=\{f(m):m\in \mathrm{Ext}_{\mathcal{M}}(c)\}.$
\end{center}

\noindent By Proposition 2.3 the graph of $f$ is also given by the formula $%
\theta (x,y,\varphi ):=\mathrm{S}(\varphi (c_{x},c_{y}))$ (where $\varphi $
is treated as a parameter). The assumption that $(\mathcal{M},S)\models
\mathrm{Repl(S)}$ then allows us to get hold of $d\in M$ such that $\mathrm{%
Ext}_{\mathcal{M}}(d)=\{f(m):m\in \mathrm{Ext}_{\mathcal{M}}(c)\}.$\medskip

\noindent \textbf{4.7.~Remark.}~An examination of the proof of Theorem B
makes it clear that the following implications hold without the assumption
of countability of $\mathcal{M}$:

\begin{center}
$(b)\Rightarrow (a)\Rightarrow (e)\Rightarrow (c)\Leftrightarrow (d).$
\end{center}

\noindent We suspect that the implication $(a)\Rightarrow (b)$ fails for
some uncountable model of ZF, but we have not been able to verify this.
However, the remaining two implications can be shown to be irreversible by
resorting to well-known uncountable models, as we shall explain.\medskip

The failure of $(e)\Rightarrow (a)$ is illustrated by the easily verified
fact that there are $\omega _{1}$-like recursively saturated models of ZF:
start with a countable recursively saturated model $\mathcal{M}_{0}$ of ZF
and let $S_{0}$ be an amenable $j$-satisfaction class $S_{0}$ for some
nonstandard $j\in \omega ^{\mathcal{M}}.$ Then use the Keisler-Morley Theorem%
\footnote{%
The Keisler-Morley Theorem is often stated for countable models of ZF, but
the usual omitting types proof of the theorem (as in \cite[Theorem 2.2.18]%
{Chang-Keisler}) works for equally well for all countable models $(\mathcal{M},S)$, where $%
\mathcal{M}$ is a model of $\mathrm{ZF}$ and $S$ is amenable over $\mathcal{M}$.} to build an $\omega _{1}$-like elementary end extension $(\mathcal{M}%
,S) $ of $(\mathcal{M}_{0},S_{0})$. It is evident that $\mathcal{M}$ is
recursively saturated but not condensable. \medskip

The failure of $(c)\Rightarrow (e)$ is illustrated by the fact that there
are $\omega _{1}$-like \textit{rather classless} recursively saturated
models of ZF; this fact was first established by Matt Kaufmann \cite%
{Kaufmann} using the combinatorial principle $\Diamond _{\omega _{1}}$;
later Shelah \cite{Shelah} used an absoluteness argument to eliminate $%
\Diamond _{\omega _{1}}$ (but no \textquotedblleft direct
proof\textquotedblright\ of this fact has been yet discovered). Note that by
Tarski's theorem on undefinability of truth, a rather classless model cannot
even carry a separative $\gamma $-satisfaction class for nonstandard $\gamma
$.$\footnote{%
Indeed by a theorem of Smith \cite{Smith}, no rather classless model $%
\mathcal{M}$ carries a $j$-satisfaction class for any nonstandard $j\in
\omega ^{\mathcal{M}}.$ Smith's result was formulated for models of PA, but
his proof works equally well for models of ZF.}$\medskip

It is also worth pointing out that the equivalence of condition (a) through
(d) of Theorem B holds for models of ZF of arbitrary cardinality provided
\textquotedblleft condensable\textquotedblright\ is replaced by
\textquotedblleft weakly condensable\textquotedblright , where $\mathcal{M}$
is said to be weakly condensable if $\mathcal{M\cong }_{p}\ \mathcal{%
M(\alpha )}\prec _{\mathbb{L}_{\mathcal{M}}}\mathcal{M}$ for some $\alpha
\in \mathrm{Ord}(\mathcal{M})$; here $\mathcal{\cong }_{p}$ denotes partial
isomorphism (two relational structures are said to be partially isomorphic
if there is a nonempty family of partial isomorphisms between them that has
the back-and-forth property, see \cite{Barwise-book}). \medskip

\noindent \textbf{4.8.~Remark.}~Condensability is a robust notion, as
indicated by (1) and (2) below.\medskip

\noindent (1) It is easy to see, using the definition of condensability,
that condensability is inherited by inner models (by an inner model of a
model $\mathcal{M}$ of ZF here we mean a transitive subclass of $\mathcal{M}$
that satisfies ZF, contains all the ordinals of $\mathcal{M}$, and is
definable in $\mathcal{M}$ by an $\mathbb{L}_{\mathcal{M}}$-formula).
\medskip

\noindent (2) The equivalence of (a) and (c) of Theorem B can be used to
show that if $\mathcal{M}$ is a condensable model of $\mathrm{ZF}$, and $%
\mathbb{P}$ is \textit{set-notion of forcing} in $\mathcal{M}$, then for
every $\mathbb{P}$-generic filter $G$ over $\mathcal{M}$, $\mathcal{M}[G]$
is also condensable (the proof is similar to the special case when $\mathcal{%
M}$ is recursively saturated, as in the proof of Theorem 2.6 of \cite%
{Enayat-Counting}). The situation is quite different for \textit{class
notions of forcing}, since as shown in Theorem 2.8 of \cite{Enayat-DO} every
countable model of ZF has a class-generic extension to a Paris model of ZF,
and of course no Paris model is condensable.

\noindent \textsc{Department of Philosophy, Linguistics, and the Theory of
Science \newline
\noindent University of Gothenburg, Gothenburg, Sweden}\newline
\noindent \texttt{email: ali.enayat@gu.se}

\end{document}